\documentstyle[11pt]{article}
\def\R{\ifmmode{\rm I\mkern-3.1mu
R\mkern1mu}\else{\rm I\kern-.18em  R\hskip1pt\ 
}\fi\relax}  

\def\inc{\subseteq} 
 
\def\b{\beta} 

\def\G{\Gamma}

\def\t{\tau}
\def\d{\delta}  
\def\th{\theta} 
\def\l{\lambda}
\def\L{\Lambda}
\def\n{\nu}
\def\di{\diamond}
\def\r{\rho}

\def\Ph{\Phi}

\def\m{\mu}
\def\s{\sigma}
\def\sou{\overline}
\def\so{\underline} 
\def\O{\Omega}
\def\f{\rightarrow}
\def\eq{\longleftrightarrow}
\def\q{\forall}
  
\def\v{\vdash}
 
\def\p{\succ}
\def\vi{\emptyset} 

\def\R{\ifmmode{\rm I\mkern-3.1mu
R\mkern1mu}\else{\rm I\kern-.18em 
R\hskip1pt\ }\fi\relax} 

\def\Z{\ifmmode{ Z\mkern-7.6mu
Z\mkern2mu}\else{ Z\kern-.28em 
Z\hskip1pt\ }\fi\relax} 

\def\Q{\ifmmode{\rm Q\mkern-10mu
l\mkern4.5mu}\else{\rm Q\kern-.57em
l\hskip3pt\ }\fi\relax} 

\def\N{\ifmmode{\rm I\mkern-3.1mu
N\mkern0.5mu}\else{\rm I\kern-.16em
N\hskip0.5pt\ }\fi\relax} 

\def\C{\ifmmode{\rm C\mkern-8.8mu
l\mkern4mu}\else{\rm C\kern-.48em
l\hskip2.6pt\ }\fi\relax} 
\def\mats{\ifmmode{ {\hbox{\bigreek s}} }\else{ 
{\bigreek s} }\fi\relax}
\def\matsin{\ifmmode{ {\hbox{\smgreek s}} }\else{ 
{\smgreek s} }\fi\relax}
\def\matt{\ifmmode{ {\hbox{\bigreek t}} }\else{ 
{\bigreek t} }\fi\relax}
\def\mattin{\ifmmode{ {\hbox{\smgreek t}} }\else{ 
{\smgreek t} }\fi\relax}

\newtheorem{theo}{Theorem}[section]
\newtheorem{lemma}{Lemma}[section]

\newtheorem{corollary}{Corollary}[section]

\begin{document} 

\Large\bf 
STORAGE OPERATORS and $\q$-POSITIVE TYPES in $TTR$ TYPE SYSTEM\\
 
\normalsize \bf 
Karim NOUR \footnote{ We thank R. David, J.L. Krivine, and M. Parigot for helpful discussions.}\\
\rm
LAMA - Equipe de Logique, Universit\'e de Chamb\'ery -
73376 Le Bourget du Lac cedex \footnote{e-mail nour@univ-savoie.fr} \\

\bf Abstract \rm 
In 1990, J.L. Krivine introduced the notion of storage operator to simulate "call by value" in the
"call by name" strategy. J.L. Krivine has shown that, using G\H{o}del translation of classical into
intuitionitic logic, we can find a simple type for the storage operators in $AF2$ type system. This
paper studies the $\q$-positive types (the universal second order quantifier appears positively in
these types), and the G\H{o}del transformations (a generalization of classical G\H{o}del translation) of
$TTR$ type system. We generalize, by using syntaxical methods, the J.L. Krivine's Theorem about
these types and for these transformations. We give a proof of this result in the case of the type
of recursive integers. \\

\bf Mathematics Subject Classification \rm: 03B40, 68Q60 \\
\bf Keywords \rm: Storage operator, Head normal form, Head reduction, $AF2$ type system, Least fixed
point, $TTR$ type system, Arrow type, Without-arrow type, $\q$-positive type, $\perp$-type, G\H{o}del
transformation.

\section{Introduction}

The strategy of left reduction (iteration of head reduction denoted by $\p$) has the following
advantages : 
\begin{itemize}
\item It has good mathematical properties stated by the normalisation Theorem : if a
$\l$-term is normalizable, then we obtain the normal form by left reduction. 
\item It seems more economic since we compute a $\l$-term only when we need it. 
\end{itemize}
Now, a drawback of the strategy of left reduction (call by name) is the fact that the argument of a
function is computed as many times as it is used. The purpose of storage operators is precisely to
correct this drawback. \\
Let $F$ be a $\l$-term (a function), and $\so{N}$ the set of normal Church integers. During the
computation, by left reduction, of $(F)\th_n$ (where $\th_n \simeq\sb{\b} \so{n}$), $\th_n$ may be
computed several times (as many times as $F$ uses it). We would like to transform $(F)\th_n$ to
$(F)\so{n}$. We also want this transformation depends only on $\th_n$ (and not $F$). In other words
we look for some closed $\l$-terms $T$ with the following properties : 
\begin{itemize}
\item For every $F$, $n \in \N$, and $\th_n \simeq\sb{\b} \so{n}$, we have $(T)\th_nF \p
(F)\so{n}$;  
\item The computation time of the head reduction $(T)\th_nF \p (F)\so{n}$  depends only
on $\th_n$.  
\end{itemize} 
Therefore the first definition : A closed $\l$-term $T$ is called storage
operator for $\so{N}$ if and only if for every $n \in \N$, and for every $\th_n \simeq\sb{\b}
\so{n}$, $(T)\th_n f \p (f)\so{n}$ (where $f$ is a new variable). \\
It is clear that a storage operator satisfies the required properties. Indeed,  
\begin{itemize}
\item Since we have $(T)\th_n f \p (f)\so{n}$, then the variable $f$ never comes in head position
during the reduction, and we may then replace $f$ by any $\l$-term. 
\item The computation time of the head reduction $(T)\th_nF \p (F)\so{n}$ depends only on $\th_n$.
\end{itemize} 
We showed (see [12]) that it is not possible to get the normal form of $\th_n$. We then change the
definition : A closed $\l$-term $T$ is called storage operator for $\so{N}$ if and only if for every 
$n \in \N$, there is a closed $\l$-term $\t_n\simeq\sb{\b} \so{n}$ (for example
$\t_n=(\so{s})^n\so{0}$, where $\so{s}$ is a $\l$-term for the successor), such that for every
$\th_n \simeq\sb{\b} \so{n}$, $(T)\th_n f \p (f)\t_n$ (where $f$ is a new variable). \\ 
If we take $T_1=\l n((n)\l x\l y(x)\l z(y)(\so{s})z)\l f(f)\so{0}$, and
$T_2=\l n\l f(((n)\l x\l y(x)(\so{s})y)f)\so{0}$, then it is easy to check that : for every $\th_n
\simeq\sb{\b} \so{n}$, $(T_1)\th_nf \p (f)(\so{s})^n\so{0}$, and $(T_2)\th_nf \p
(f)(\so{s})^n\so{0}$. Therefore $T_1$ and $T_2$ are storage operators for $\so{N}$.\\

The $AF2$ type system is a way of interpreting the proof rules for the second order intuitionistic
logic plus equational reasoning as construction rules for terms. In this system we have the
possibility to define the data types, the representation in $\l$-calculus being automaticaly extracted
from the logical definition of the data type. At the logical level the data type are defined by
second order formulas expressing the usual iterative definition of the corresponding algebras of
terms and the data receive the corresponding iterative definition in $\l$-calulus. For example, the
type of integers is the formula : $N[x]=\q X \{\q y[X(y) \f X(sy)] \f [X(0) \f X(x)] \}$ ($X$ is a
unary predicate variable, $0$ is a constant symbol for zero, and $s$ is a unary function symbol for
successor). \\
If we try to type a storage operator $T$ in $AF2$ type system, we naturally find the type
$\q x \{ N[x] \f [(N[x] \f O)\f O] \}$ (where $O$ is a particular $0$-ary predicate symbol
which represents an arbitrary type). Indeed, if $\v_{AF2} \t_n:N[s^n(0)]$, and $f$ is of type
$N[s^n(0)] \f O$, then $f:N[s^n(0)] \f O \v_{AF2} (f)\t_n:O$. It is natural to have
$(T)\th_nf$ of type $O$. If $\v_{AF2}\th_n:N[s^n(0)]$, then the type for $T$ must be
$\q x \{ N[x] \f [(N[x]\f O)\f O] \}$.\\
It is easy to check that $\v_{AF2}T_1,T_2:\q x \{ N[x] \f [(N[x]\f O)\f O] \}$. \\
The type $\q x \{ N[x] \f [(N[x]\f O)\f O] \}$ does not characterize the storage operators.
Indeed, if we take $T=\l n\l f(f)n$, we obtain : 
\begin{itemize}
\item $n:N[x],f:N[x]\f O \v_{AF2}(f)n:O$, then, $\v_{AF2}T:\q x \{ N[x] \f [(N[x] \f O)\f
O] \}$.  
\item For every $\th_n \simeq\sb{\b} \so{n}$, $(T)\th_nf \p (f)\th_n$, therefore $T$ is not a storage
operator for $\so{N}$.  
\end{itemize}
This comes from the fact that the type $\q x \{ N[x] \f [(N[x]\f O)\f O] \}$ does not take
into account the independance of $\t_n$ with $\th_n$. To solve this problem, we must prevent the use
of the first $N[x]$ in $\q x \{ N[x] \f [(N[x]\f O)\f O] \}$ as well as his subtypes to
prove the second. Therefore, we will replace the first $N[x]$ by a new type $N$*$[x]$ with the
following properties :
\begin{itemize}
\item  $\v_{AF2}\so{n}:N$*$[s^n(0)]$ (for example, take $N$*$[x]=\q X \{ \q y[F(X,y) \f F(X,sy)] \f
[F(X,0) \f F(X,x)] \}$) ; 
\itemÊIf $\n:N$*$[x], x_i:\q y[F(G,y) \f F(G,sy)], y_j:F(H,a) \v_{AF2} t:N[s^n(0)]$, then
$\v_{AF2}t':N[s^n(0)]$, where $t'$ is the normal form of $t$ ; 
\item There is a closed $\l$-term $T$, such that $\v_{AF2}T:\q x \{N$*$[x] \f [(N[x] \f O)\f
O] \}$.
\end{itemize} 
A simple solution for the second property is to take a formula $F(X,a)$ ending with a new constant
symbol. Indeed, since $N[x]$ does not contain this symbol, we cannot use the variables $\n,x_i,y_j$
in the typing of $t'$. We suggest the following proposition : 
\begin{center}
$N$*$[x]=\q X \{\q y[(X(y) \f O) \f (X(sy) \f O)] \f [(X(0) \f O) \f (X(x) \f O)] \}$.
\end{center}
It is easy to chech that $\v_{AF2}T_1,T_2:\q x \{N$*$[x] \f [(N[x]\f O) \f O] \}$ (see [6] and
[12]).\\ 
For each formula $F$ of $AF2$, we indicate by $F$* the formula obtained by putting $\neg$ in front of
each atomic formulas of $F$ ($F$* is called the G\H{o}del translation of $F$).\\ 
J.L. Krivine has shown  that the type $\q x \{ N$*$[x] \f \neg \neg N[x] \}$ characterize the storage
operators for $\so{N}$ (see [6]). But the $\l$-term $\t_n$ obtained may contain variables substituted
by $\l$-terms $u_1,...,u_m$ depending on $\th_n$. Since the $\l$-term $\t_n$ is $\b \eta$-equivalent
to $\so{n}$, therefore, the left reduction of the $\t_n[u_1/x_1,...,u_m/x_m]$ is equivalent to the
left reduction of $\t_n$ and the $\l$-terms $u_1,...,u_m$ will therefore never be evaluated during
the reduction. \\ Taking into account the above remarks, we modify again the definition : A closed
$\l$-term $T$ is called a storage operator for $\so{N}$ if and only if for every $n \in \N$, there is
a $\l$-term $\t_n\simeq\sb{\b} \so{n}$, such that for every $\th_n \simeq\sb{\b} \so{n}$, there is a
substitution $\s$, such that $(T)\th_nf \p (f)\s(\t_n)$ (where $f$ is a new variable). \\

The $AF2$ type system is satisfactory from an extensional point of view : one can construct programs
for all the functions whose termination is provable in the second order Peano arithmetic. But from an
intensional point of view the situation is very different : we cannot always obtain the simple (in
term of time complexity, for instance) programs we need. For example we cannot find a $\l$-term of
type $\q x \q y \{ N[x],N[y] \f N[min(x,y)] \}$ ($min$ is a binary function symbol defined by
equations) in $AF2$ type system that computes the minimum of two Church integers in time $O(min)$
\footnote{ R. David gives a $\l$-term of type $N,N \f N$ $(N=\q X \{ [X \f X] \f [X \f X] \})$ in $F$
type system that computes the minimum of two Church integers in time $O(min.Log(min))$. The notion of
storage operators plays an important tool in this constraction (see [2]).} \\ 
The $TTR$ type system is an extension of $AF2$ based on recursive definitions of types, which is
intented to solve the basic problems of efficiency mentioned before. In $TTR$ we have a logical
operator $\m$ of least fixed point. If $A$ is a formula, $C$ an $n$-ary predicate symbol which
appears and occurs positively in $A$, $x_1,...,x_n$ first order variables, and $t_1,...,t_n$ terms,
then $\m C x_1...x_n A<t_1,..,t_n>$ is a formula called the least fixed point of $A$ in $C$
calculated over the terms $t_1,...,t_n$. The interded logical meaning of the formula $\m C x_1...x_n
A<t_1,..,t_n>$ is $K(t_1,...,t_n)$, where $K$ is the least $X$, such that $X(x_1,...,x_n) \eq A$.
$TTR$ allows to define the multisorted term algebras as least fixed points. For example the type of
recursive integers is the formula : $N^r[x]=\m Cz[\q X \{ \q y[C(y) \f X(sy)] \f [X(0) \f X(z)]
\}]<x>$ ($X$ is a unary predicate variable, $0$ is a constant symbol for zero, and $s$ is a unary
function symbol for successor). \\

In this paper we study the types $D$ of $TTR$, and the transformations *, for which we have the
following result : if $\v_{TTR}T:D$*$ \f \neg \neg D$, then for every $\l$-term $t$ with
$\v_{TTR}t:D$, there are $\l$-terms $\t_t$ and $\t'_t$ such that $\t_t \simeq\sb{\b} \t'_t$,
$\v_{TTR}\t'_t:D$, and for every $\th_t \simeq\sb{\b} t$, there is a substitution $\s$, such that
$(T)\th_tf \p (f)\s(\t_t)$ (where $f$ is a new variable). \\ 
We prove \footnote{J.L. Krivine and the author proved independely the same result for $AF2$ type
system (see [7] and [12]).} that, to obtain this result, it suffies to assume that :
\begin{itemize}
\item  The universal second order quantifier appears positively in $D$ ($\q$-positive type)
\footnote{ This types were studied by some authors (in particular R. Labib-Sami), and have
remarkable properties (see [8]).}.  
\item  The transformation * satisfies the following properties : 
\begin{itemize}
\item[]Ê- If $A=C(t_1,...,t_n)$, then $A$*$=A$ ; 
\item[]Ê- If $A=X(t_1,...,t_n)$, then $A$*$=F_X[t_1/x_1,...,t_n/x_n]<X_1,...,X_r>$ where $F_X$ is a
formula ending with $\perp$ and having $x_1,...,x_n,X_1,...,X_r$ as free variables ;  
\item[]Ê- $(A \f B)$*$=A$*$ \f B$* ;
\item[]Ê- $(\q xA)$*$=\q xA$*.  
\item[]Ê- $(\q XA)$*$= \q X_1...X_rA$*. 
\item[]Ê- $(\m C x_1...x_n A<t_1,..,t_n>)$*$=\m C x_1...x_n A$*$<t_1,..,t_n>$.
\end{itemize}
\end{itemize}

We give the proof of this result in the case of the type of recurcive integers.

\section{ Basic notions of pure $\l$-calculus}

Our notation is standard (see [1] and [5]). \\
We denote by $\L$ the set of terms of pure $\l$-calculus, also called $\l$-terms.\\
Let $t,u,u_1,...,u_n \in \L$, the application of $t$ to $u$ is denoted by $(t)u$. In the same way we
write $(t)u_1...u_n$ instead of $(...((t)u_1)...)u_n$. \\
The $\b$-reduction (resp. $\b$-equivalence) is denoted by $t \f_{\b} u$ (resp. $t \simeq\sb{\b} u$).\\
The set of free variables of a $\l$-term $t$ is denoted by $Fv(t)$. \\
The notation $t[u_1/x_1,...,u_n/x_n]$ represents the result of the simultaneous substitution of
$\l$-terms $u_1,...,u_n$ to the free variables $x_1,...,x_n$ of $t$ (after a suitable renaming of the
bounded variables of $t$). \\
With each normal $\l$-term, we associate a set of $\l$-terms $STE(t)$ by induction : \\
if $t=\l x_1...\l x_n(y)t_1...t_m$, then $STE(t)=\{ t \} \bigcup \displaystyle\bigcup_{1 \leq i \leq
m}STE(t_i)$.\\  
Let us recall that a $\l$-term $t$ either has a head redex [i.e. $t=\l x_1...\l
x_n(\l xu)vv_1 ... v_m$, the head redex being $(\l xu)v$], or is in head normal form [i.e. $t=\l
x_1...\l x_n(x)vv_1... v_m$]. \\ 
The notation $t \p t'$ means that $t'$ is obtained from $t$ by some head reductions, and we denote by
$n(t,t')$, the number of steps to go from $t$ to $t'$. \\
A $\l$-term $t$ is said to be solvable if and only if the head reduction of $t$ terminates. \\
We define an equivalence relation $\sim$ on $\L$ by : $u \sim v$ if and only if there is a $t$, such
that $u \p t$, and $v \p t$. In particular, if $v$ is in head normal form, then $u \sim v$ means that
$v$ is the head normal form of $u$. 

\begin{theo} ([6]). If $t \p t'$, then for every $u_1,...,u_n \in \L$ :\\
1) there is a $v \in \L$, such that $(t)u_1...u_n \p v$, $(t')u_1...u_n \p v$, and
$n((t)u_1...u_n,v)= n((t')u_1....u_n,v)+n(t,t')$.\\ 
2) $t[u_1/x_1,...,u_n/x_n] \p t'[u_1/x_1,...,u_n/x_n]$, and
$n(t[u_1/x_1,...,u_n/x_n],t'[u_1/x_1,...,u_n/x_n])= n(t,t')$. 
\end{theo}

\bf Remark. \rm Theorem 2.1 shows that to make the head reduction of $(t)u_1...u_n$ (resp.
$t[u_1/x_1,...,u_n/x_n]$), it is equivalent (same result, and same number of steps) to make some steps
in the head reduction of $t$, and then make the head reduction of $(t')u_1...u_n$ (resp.
$t'[u_1/x_1,...,u_n/x_n]$).

\section{Basic notions of typed $\l$-calculus}

\subsection{The $AF2$ type system}

The types will be formulas of second order predicate logic over a given language. \\
The logical symbols are $\perp$ (for absurd), $\f$ and $\q$ (and no other ones).\\ 
There are individual variables : $x,y,...$ (also called first order variables) and $n$-ary predicate
variables ($n=0,1,...$) : $X,Y,...$ (also called second order variables).\\ 
The terms and the formulas are up in the usual way.\\ 
The formula $F_1 \f (F_2 \f(...\f(F_n \f G)...))$ is denoted by $F_1,F_2,...,F_n \f G$, and $F \f
\perp$ is denoted by $\neg F$. The formula $\q v_1...\q v_n F$ is denoted by $\q \bf v$\rm$F$, and
the sentence "\bf v \rm is not free in $A$" means that for all $1 \leq i \leq n$, $v_i$ is not free in
$A$. \\
If $X$ is a unary predicate variable, $t$ and $t'$ two terms, then the formula $\q X[Xt \f Xt']$ is
denoted by $t=t'$, and is said to be equation. A particular case of $t=t'$ is a formula of the forme
$t[u_1/x_1,...,u_n/x_n]=t'[u_1/x_1,...,u_n/x_n]$ or $t'[u_1/x_1,...,u_n/x_n]=t[u_1/x_1,...,u_n/x_n]$,
$u_1,...,u_n$ being terms of the language.  \\
After, we denote by \bf E \rm a system of function equations.\\ 
A context $\G$ is a set of the form $x_1:A_1,...,x_n:A_n$ where $x_1,...,x_n$ are distinct variables
and $A_1,...,A_n$ are formulas.  \\
We are going to describe a system of typed $\l$-calculus called second order functional arithmetic
(shortened in $AF2$ for Arithm\'etique Fonctionnelle du seconde ordre). The rules of typing are the
following :  \begin{itemize}
\item []Ê(1) $\G,x:A \v_{AF2} x:A$. 
\item [] (2) If $\G,x:B \v_{AF2} t:C$, then $\G \v_{AF2} \l xt:B \f C$. 
\item [] (3) If $\G \v_{AF2}u:B \f C$, and $\G \v_{AF2}v:B$, then $\G \v_{AF2}(u)v:C$.
\item [] (4) If $\G \v_{AF2}t:A$, and $x$ does not appear in $\G$, then $\G \v_{AF2}t:\q xA$. 
\item [] (5) If $\G \v_{AF2}t:\q xA$, then, for every term $u$, $\G \v_{AF2}t:A[u/x]$. 
\item [] (6) If $\G \v_{AF2}t:A$, and $X$ does not appear in $\G$, then $\G \v_{AF2}t:\q XA$. 
\item [] (7) If $\G \v_{AF2}t:\q XA$, then, for every formula $G$, $\G
\v_{AF2}t:A[G/X(x_1,...,x_n)]$ (*)  
\item [] (8) If $\G \v_{AF2}t:A[u/x]$, then $\G \v_{AF2}t:A[v/x]$, $u=v$ being a particular case of an
equation of \bf E \rm. 
\end{itemize}

(*) $A[G/X(x_1,...,x_n)]$ is obtained by replacing in $A$ each atomic formula $X(t_1,...,t_n)$ by\\ 
$G[t_1/x_1,...,t_n/x_n]$. To simplify, we write sometimes $A[G/X]$ instead of
$A[G/X(x_1,...,x_n)]$.\\

Whenever we obtain the typing $\G \v_{AF2}t:A$ by means of these rules, we say that "the $\l$-term
$t$ is of type $A$ in the context $\G$, with respect to the equation of \bf E \rm". 

\begin{theo} ([5],[9]).\\
1)  Conservation Theorem: If $\G \v_{AF2}t:A$, and $t \f_{\b} t'$, then $\G \v_{AF2}t':A$.\\
2)  Strong normalization: If $\G \v_{AF2}t:A$, then $t$ is strongly normalizable.
\end{theo}

\subsection{The $TTR$ type system}

Let $X$ be a predicate variable or predicate symbol, and $A$ a type of $AF2$. \\
We define the notions "$X$ is positive in $A$" and "$X$ is negative in $A$" by induction :
\begin{itemize}
\item[] - If $X$ does not appears in $A$, then $X$ is positive and negative in $A$ ;
\item[] - If $A=X(t_1,...,t_n)$, then $X$ is positive in $A$, and $X$ is not negative in $A$ ;
\item[] - If $A=B \f C$, then $X$ is positive (resp. negative) in $A$ if and only if $X$ is negative
(resp. positive) in $B$, and $X$ is positive (resp. negative) in $C$ ; 
\item[] - If $A=\q vB$, and $v \not = X$, then $X$ is positive (resp. negative) in $A$ if and only if
$X$ is positive (resp. negative) in $B$.
\end{itemize} 

We add to the second order predicate calculus a new logic symbol $\m$, and we allow a new
construction for formulas : if $A$ is a formula, $C$ an $n$-ary predicate symbol which appears
positively in $A$, $x_1,...,x_n$ first order variables, and $t_1,...,t_n$ terms, then
$\m C x_1...x_n A<t_1,...,t_n>$ is a formula called the least fixed point of $A$ in $C$ calculated
over the terms $t_1,...,t_n$. \\
We extend the notions "$X$ is positive in a type" and "$X$ is negative in a type" by the following way
: $X$ is positive (resp. negative) in $\m C x_1...x_n A<t_1,...,t_n>$ if and only if $X$ is positive
(resp. negative) in $A$. \\
We extend the definition of the substitution by assuming that $C,x_1,...,x_n$ are
bounded in the formula $\m C x_1...x_n A<t_1,...,t_n>$.\\
We define on these formulas a binary relation $\inc$ by : $A \inc B$ if and only if it is obtained by
using the following rules :\\

\begin{minipage}[t]{250pt}
$  (ax)  A \inc A$\\
\end{minipage} 
\begin{minipage}[t]{250pt}\small\sl 
$  (\f )   \displaystyle\frac{ A \inc A' \quad B \inc B' } { A' \f B \inc A \f B' }$
\\ \end{minipage}

\begin{minipage}[t]{250pt}
$  (\q i_g)  \displaystyle\frac{ A[G/v] \inc B } { \q vA \inc B }$ (1)\\
\end{minipage} 
\begin{minipage}[t]{250pt}\small\sl 
$  (\q i_d)   \displaystyle\frac{ A \inc B } { A \inc \q vB }$ (2) \\
\end{minipage}

\begin{minipage}[t]{250pt}
$  ( e )  \displaystyle\frac{ A \inc B[v/y] } { A \inc B[w/y] }$ (3)\\
\end{minipage} 
\begin{minipage}[t]{250pt}\small\sl 
$  ( tr)   \displaystyle\frac{ A \inc D \quad D \inc B } {A \inc B }$ \\
\end{minipage}

$(\m_d)$   $D[\m C x_1...x_m D<z_1,...,z_m> / C(z_1,...,z_m)][t_1/x_1,...,t_m/x_m] \inc \m C x_1...x_m
D<t_1,...,t_m>$  \\

$(\m'_g)$  $\m C x_1...x_m D<t_1,...,t_m> \inc D[\m C x_1...x_m D<z_1,...,z_m> /
C(z_1,...,z_m)][t_1/x_1,...,t_m/x_m]$

\begin{center}
$(\m_g)$   $\displaystyle\frac{D[E/C(x_1,...,x_m)] \inc E} {\m C x_1...x_m
D<t_1,...,t_m> \inc E[t_1/x_1,...,t_m/x_m]}$ 
\end{center}

(1) $G$ is a formula if $v$ is a second order variable, and a term if $v$ is a first order
variable.\\  
(2) $v$ is not free in A.\\ 
(3) $v=w$ is a particular case of an equation of \bf E \rm.\\

$(\m_d)$ and $(\m'_g)$ are the rules of factorisation and development of a fixed point.\\
$(\m_g)$ expresses the fact that $\m C x_1...x_m D<t_1,...,t_m>$ is a least fixed point.\\

We are going to describe a system of typed $\l$-calculus called theory of recursive types (shortened
in $TTR$ for Th\'eorie des Types R\'ecursifs) where the types are formulas of language. The rules of
typing are the following :
\begin{itemize} 
\item[] - The typing rules (1),...,(8) of $AF2$ type system.
\item[] - $(\inc )$  $\displaystyle\frac{\G \v_{TTR}t:A \quad A \inc B} {\G \v_{TTR}t:B}$ 
\item[] - $(Y)$  $\displaystyle\frac{\G \v_{TTR}t:\q x_1...\q x_m [C(x_1,...,x_m) \f E] \f \q x_1...\q
x_m [D \f E] }{\G \v_{TTR}(Y)t:\q x_1...\q x_m[\m C x_1...x_m D<x_1,...,x_m> \f E]}$  \\ 	 
where $C$ is not free in $E$ and $G$, and $Y$ is the Turing's fixed point. 
\end{itemize}
The rule $(Y)$  expresses also the fact that $\m C x_1...x_m D<t_1,...,t_m>$ is a least fixed point.

\begin{theo} ([12],[18]).\\ 
1) Conservation Theorem If $\G \v_{TTR}t:A$, and $t \f_{\b} t'$, then $\G \v_{TTR}t':A$.\\
2) Strong normalization If $\G \v_{TTR}t:A$ without using the rule $(Y)$, then $t$ is strongly
normalizable.\\ 
3) Weak normalization If $\G \v_{TTR}t:A$, and if all least fixed points of $A$ are
positives, then $t$ is normalizable.
\end{theo}

The $TTR^{\di}$ type system is the subsystem of $TTR$ where we only have propositional variables and
constants (predicate variables or predicate symbols are of arity 0). So, first order variables,
function symbols, and finite sets of equations are useless. With each predicate variable (resp.
predicate symbol) $X$, we associate a predicate variable (resp. a predicate symbol) $X^{\di}$ of
$TTR^{\di}$ type system. For every formula $A$ of $TTR$, we define the formula $A^{\di}$ of
$TTR^{\di}$ obtained by forgetting in $A$ the first order part. If $\G =x_1:A_1,...,x_n:A_n$ is a
context of $TTR$, then we denote by $\G^{\di}$, the context $x_1:A_1^{\di},...,x_n:A_n^{\di}$ of
$TTR^{\di}$. We write $\G \v_{TTR^{\di}}t:A$ if $t$ is tyable in $TTR^{\di}$ of type $A$ in the
context $\G$.

\begin{theo} If $\G \v_{TTR}t:A$, then $\G^{\di} \v_{TTR^{\di}}t:A^{\di}$. 
\end{theo} 
\bf Proof \rm By induction on the length of the derivation $\G \v_{TTR}t:A$. $\Box$

\begin{theo} \quad \\
1) Conservation Theorem If $\G \v_{TTR^{\di}}t:A$, and $t \f_{\b} t'$, then $\G \v_{TTR^{\di}}t':A$.\\
2) Strong normalization If $\G \v_{TTR^{\di}}t:A$ without using the rule $(Y)$, then $t$ is strongly
normalizable.\\
3) Weak normalization If $\G \v_{TTR^{\di}}t:A$, and if all least fixed points of $A$ are
positives, then $t$ is normalizable.
\end{theo} 
\bf Proof \rm We use Theorems 3.2 and 3.3. $\Box$ \\

\bf Remark \rm We cannot if the reverse of 2)-Theorem 3.2 is true, but the $\l$-term \\
$t=\l x(\l y((x)(y)\l xx)(y)\l x\l yx)\l x(x)x$ (which is strongly normalizable, and untypable in
$AF2$ type system (see [3])) is typable in $TTR$ type system. Indeed, if we take $B=\m C(\q XX \f
C)$, we check easily that $\v_{TTR^{\di}}t:[B \f (B \f B)] \f B$.

\section{ Properties of $TTR$ type system}

\subsection{ Permutations Lemmas }

\begin{lemma} 1) The typing rules (5), (7), and (8) are admissible.\\
2) In the typing, we may replace the succession of $n$ times ($\inc$) and $m$ times (4) (resp. (6)),
by the succession of $m$ times (4) (resp. (6)) and n times ($\inc$). \\
3) If $\G \v_{TTR}t:B$ is derived from $\G \v_{TTR}t:A$, then we may assume that we begin by the
applications of (4), (6), and next by ($\inc$). 
\end{lemma}  
\bf Proof \rm Easy. $\Box$

\begin{lemma} 1) If $A \inc B$, then, for every sequence of terms and/or formulas \bf G, \it
$A \bf [G/v] \it \inc B\bf[G/v]$, \it and we use the same proof rules.\\ 
2) If $\G \v_{TTR}t:A$, then, for every sequence of terms and/or formulas \bf G, \it
$\G \bf[G/v] \it \v_{TTR}t:A \bf [G/v]$, \it and we use the same typing rules. 
\end{lemma}  
\bf Proof \rm  By induction on the length of the derivation $A \inc B$ (resp. $\G \v_{TTR}t:A$).
$\Box$

\begin{corollary} If $\G,x:A \v_{TTR}(x)u_1...u_n:B$, then :\\
$n=0$, and there is $\bf v_0$ \it not free in $A$ and $\G$, such that $\q \bf v_0$\it$A \inc B$,\\
or \\
$n \geq 1$, and there are types $C_i$,$B_i$ $(i=1,...,n)$ and $\bf v_i$\it$(i=1,Én)$ not free in $A$
and $\G$, such that $\q \bf v_0$\it$A \inc C_1 \f B_1$, $ \q \bf v_i$\it$B_i \inc C_{i+1} \f B_{i+1}$
$1 \leq i \leq n-1$, $\q \bf v_n$\it$B_n \inc B$, and $\G,x:A \v_{TTR} u_i:C_i$ $1\leq i \leq n$.
\end{corollary}  
\bf Proof \rm  By induction on $n$. $\Box$

\begin{lemma} 1) If $X$ is positive (resp. negative) in $D$, and $A \inc B$, then $D[A/X] \inc
D[B/X]$ (resp. $D[B/X] \inc D[A/X]$).\\ 
2) We may eliminate the rule ($\m'_g$).
\end{lemma} 
\bf Proof \rm 1) By induction on $D$. \\
2) By rule ($\m_d$), we have $A[\m C x_1...x_n A<y_1,...,y_n>/C(y_1,...,y_n)] \inc \m C x_1...x_n
A<x_1,...,x_n>$, then, by 1), $A[A[\m C x_1...x_n A<y_1,...,y_n>/C(y_1,...,y_n)]/C(x_1,...,x_n)] \inc
A[\m C x_1...x_n A<x_1,...,x_n>/C(x_1,...,x_n)]$,  and, by using the rule ($\m_g$), we obtain
$\m C x_1...x_n A<t_1,...,t_n> \inc A[\m C x_1...x_n
A<y_1,...,y_n>/C(y_1,...,y_n)][t_1/x_1,...,t_n/x_n]$. $\Box$

\subsection{Without-arrow types and arrow types} 

\bf Definitions \rm \\
1) A type $A$ is said to be without-arrow type if and only if $A$ does not contain
any arrow. \\
2) Each without-arrow type $A$ contains a unique atomic formula $X(t_1,...,t_n)$. We denote $X$ by
$At(A)$. We distinguish between two kinds of without-arrow types :\\ 
- A without-arrow type $A$ is said to be of kind 1 if and only if $At(A)$ is free in $A$.\\ 
- A without-arrow type $A$ is said to be of kind 2 if and only if $At(A)$ is bounded in $A$.

\begin{lemma} 1) If $A$ is a without-arrow type of kind 1, and $A \inc B$, then $B$ is a without-arrow
type of kind 1, and $At(A)=At(B)$.\\
2) If $A$ is a without-arrow type of kind 2, then, for every type $B$, we have $A \inc B$.
\end{lemma}
\bf Proof \rm 1) By induction on the length of the derivation $A \inc B$.\\ 
2) Easy. $\Box$ \\

\bf Definition \rm  A type $A$ is said to be arrow type if and only if $A$ contains at least an
arrow.

\begin{lemma} If $A$ is an arrow type, and $A \inc B$, then $B$ is an arrow type.
\end{lemma}
\bf Proof \rm  By induction on the length of the derivation $A \inc B$. $\Box$

\begin{corollary} Let $A$ be an atomic formula. If $\G \v_{TTR}t:A$, then $t$ does not begin by $\l$.
Other words, if $\G \v_{TTR}\l xu:B$, then $B$ is an arrow type. 
\end{corollary}
\bf Proof \rm If $t$ begins by $\l$, then there are $E,F$, and \bf v, \rm such that $\q \bf
v$\rm$(E \f F) \inc A$, therefore, by Lemma 4.5, $A$ is an arrow type. $\Box$ \\
                                     
\bf Definition \rm  For every arrow type $A$, we define the type $Rep(A)$ as follows, by induction on
A : 
\begin{itemize}
\item[] - $Rep(E \f F)=E \f F$ ;  
\item[] - $Rep(\q vB)=\q vRep(B)$ ;
\item[] - $Rep(\m C x_1...x_n B<t_1,...,t_n>)= \\
 Rep(B)[\m C x_1...x_n B<y_1,...,y_n>/C(y_1,..,y_n)][t_1/x_1,...,t_n/x_n]$. 
\end{itemize}

\begin{lemma} If $A$ is an arrow type, then : \\
1) there are $G,D$ and \bf v \it such that $Rep(A)=\q \bf v \it (G \f D)$. \\
2) $A \inc Rep(A)$, and $Rep(A) \inc A$.
\end{lemma}
\bf Proof \rm By induction on $A$. $\Box$ \\

\bf Remark. \rm The Lemma 4.6 means that if $A$ is an arrow type, then $Rep(A)$ is an "equivalent"
type to $A$ of the form $\q \bf v$\rm$(G \f D)$. In the rest of the paper, we denoted $G$ by $A_g$ and
$D$ by $A_d$.

\begin{lemma} Let $A,B$ be two types, and $X,X'$ two predicate variables or predicate symbols, such
that $X'$ is not free in $A$.\\
1) If $X$ is positive in $A$, and $X'$ is positive in $B$, then $X'$ is positive in $A[B/X]$. \\ 
2) If $X$ is positive in $A$, and $X'$ is negative in $B$, then $X'$ is negative in $A[B/X]$. \\
3) If $X$ is negative in $A$, and $X'$ is positive in $B$, then $X'$ is negative in $A[B/X]$. \\
4) If $X$ is negative in $A$, and $X'$ is negative in $B$, then $X'$ is positive in $A[B/X]$. 
\end{lemma}
\bf Proof \rm By induction on $A$. $\Box$

\begin{lemma} Let $A$ be an arrow type. \\
1) If $X$ is positive (resp. negative) in $A$, then $X$ is positive (resp. negative) in $Rep(A)$. \\
2) If \bf G \it is a sequence of terms and/or formulas, then $Rep(A \bf [G/v])\it= Rep(A) \bf
[G/v]$. 
\end{lemma}
\bf Proof \rm 1) We argue by induction on $A$. The only non-trivial case is the one where
$A=\m C x_1...x_n B<t_1,...,t_n>$. If $X$ is positive (resp. negative) in $A$, then $X$ is positive
(resp. negative) in $B$. By the induction hypothesis, we have $X$ is positive (resp. negative) in
$Rep(B)$, therefore, by Lemma 4.7, $X$ is positive (resp. negative) in $Rep(A)$. \\     
2) By induction on $A$. $\Box$
                                                                                                     
\begin{theo} Let $A,B$ be two arrow types, such that $Rep(A)=\q \bf v \it (A_g \f A_d)$ and $Rep(B)=
\q \bf v' \it (B_g \f B_d)$. If $A \inc B$, then there is a sequence of terms and/or formulas \bf G,
\it such that $B_g \inc A_g \bf [G/v],$ \it and $A_d \bf [G/v] \it \inc B_d$. 
\end{theo}
\bf Proof \rm  We argue by induction on the length of the derivation $A \inc B$. Let us look at the
rule used in the last step. The only non-trivial cases are : 
\begin{itemize}                                                                                                                                   
\item[] -  (tr) : then $A \inc D$, and $D \inc B$. If $Rep(D)=\q \bf v"$\rm$(D_g \f D_d)$, by the
induction hypothesis, there are sequences \bf G \rm and \bf G$"$ \rm such that $D_g \inc A_g \bf
[G/v]$, \rm $A_d \bf [G/v]$\rm$\inc D_d$,  $B_g \inc D_g \bf [G"/v"]$, \rm and $D_d \bf [G"/v"]$\rm$
\inc B_d$. It is clear that we may assume that \bf v$"$ \rm is not free in $A_g$ and $A_d$, therefore,
by Lemma 4.2, we have $B_g \inc A_g \bf [G/v][G"/v"]$, \rm and $A_d \bf [G/v][G"/v"]$\rm$\inc B_d$.
Let $\bf G'=G[G"/v"]$, \rm then $B_g \inc A_g \bf [G'/v]$, \rm and $A_d \bf [ G'/v]$\rm$\inc B_d$.   
\item[] -  $(\m_d)$ : then $A=D[\m C x_1...x_k D<y_1,...,y_k>/C(y_1,...,y_k)][t_1/x_1,..,t_k/x_k]$,
and $B=\m C x_1...x_k D<t_1,...,t_k>$. Therefore, by Lemma 4.8, $Rep(A)=Rep(B)$, $A_g=B_g$, and
$B_d=A_d$, and so $B_g \inc A_g$, and $A_d \inc B_d$.  
\item[] -  $(\m_g)$ : then $A=\m C x_1...x_k D<t_1,...,t_k>$, $B=E[t_1/x_1,...,t_k/x_k]$, and
$D[E/C(x_1,...,x_k)] \inc E$. Therefore $Rep(D)=\q \bf v$\rm$(D_g \f D_d)$ with \\
$D_g[\m C x_1...x_k D<y_1,...,y_k>/C(y_1,...,y_k)][t_1/x_1,..,t_k/x_k]=A_g$, \\
$D_d[\m C x_1...x_k D<y_1,...,y_k>/C(y_1,...,y_k)][t_1/x_1,..,t_k/x_k]=A_d$, and  \\
$Rep(E)=\q \bf v'$\rm$(E_g \f E_d)$ with $E_g[t_1/x_1,...,t_k/x_k]=Bg$,
$E_d[t_1/x_1,...,t_k/x_k]=B_d$.  \\
By the induction hypothesis, there is a sequence \bf G, \rm such that $E_g \inc D_g[E/C(x_1,...,x_k)]
\bf [G/v]$ \rm, and $D_d[E/C(x_1,...,x_k)] \bf [G/v]$\rm$\inc E_d$. $C$ is positive in $D$, therefore,
by Lemma 4.8, $C$ is negative in $D_g$, and $C$ is positive in $D_d$.\\  
$D[E/C(x_1,...,x_k)] \inc E$, then $\m C x_1...x_k D<y_1,...,y_k> \inc E[y_1/x_1,...,y_k/x_k]$, and,
by 1)-Lemma 4.3, $Eg \inc Dg[\m C x_1...x_k D<y_1,...,y_k>/C(y_1,...,y_k)] \bf [G/v]$ \rm, and \\  
$D_d[\m C x_1...x_k D<y_1,...,y_k>/C(y_1,...,y_k)]\bf [G/v]$\rm$\inc E_d$, and so, by Lemma 4.2,\\
$E_g[t_1/x_1,...,t_k/x_k] \inc D_g[\m C x_1...x_k D<y_1,...,y_k>/C(y_1,...,y_k)]\bf
[G/v]$\rm$[t_1/x_1,...,t_k/x_k]$, and $D_d[\m C x_1...x_k D<y_1,...,y_k>/C(y_1,...,y_k)] \bf
[G/v]$\rm$[t_1/x_1,...,t_n/x_n] \inc E_d[t_1/x_1,...,t_k/x_k]$. Let $\bf G'=G$\rm$
[t_1/x_1,...,t_k/x_k]$, then $B_g \inc A_g \bf [G'/v]$,\it and $A_d \bf [G'/v]$\rm$\inc B_d$. $\Box$ 
\end{itemize}
                                                                                
\begin{corollary} Let $B$ be an atomic formula. If $\G,x:A \f B \v_{TTR}(x)u_1...u_n:C$, then $n=1$,
and $\G,x:A \f B \v_{TTR}u_1:A$.
\end{corollary}
\bf Proof \rm By Corollary 4.1, we have $\q \bf v$\rm$(A \f B) \inc F \f G$, $\G,x:A \f
B \v_{TTR}u_1:F$, and \bf v \rm is not free in $\G$ and $A \f B$. Therefore, by Theorem 4.1, $F \inc
A$, and $B \inc G$, then $\G,x:A \f B \v_{TTR}u_1:A$. If $n>1$, then $\q \bf v'$\rm$G \inc H \f J$,
and $\bf v'$ \rm is not free in $\G$ and $A \f B$. Therefore $\q \bf v'$\rm$B \inc H \f J$, and $\q
\bf v'$\rm$B$ is a without-arrow type of kind 1. A contradiction. $\Box$ 

\begin{lemma} If $x_1:A_1,...,x_n:A_n \v_{TTR}t:A$, $B_i \inc A_i$ $1 \leq i \leq n$, and $A \inc B$,
then $x_1:B_1,...,x_n:B_n \v_{TTR}t:B$. 
\end{lemma}
\bf Proof \rm We argue by induction on $t$. The only non-trivial cases are :
\begin{itemize}
\item[]Ê- If $t=\l xu$, then $x_1:A_1,...,x_n:A_n,x:E \v_{TTR}u:F$, $\q \bf v$\rm$(E \f F) \inc A$,
and \bf v \rm is not free in $E$ and $A_j$ $1 \leq j \leq n$. We may assume that \bf v \rm is not
free in $E$ and $B_j$ $1 \leq j \leq n$. By the induction hypothesis, we have
$x_1:B_1,...,x_n:B_n,x:E \v_{TTR}u:F$, and so $x_1:B_1,...,x_n:B_n \v_{TTR}t:B$.  
\item[] - If $t=(Y)u$, then $\q \bf v$\rm$\q y_1...\q y_m[\m C y_1...y_m E<y_1,...,y_m> \f D]) \inc
A$,  $x_1:A_1,...,x_n:A_n \v_{TTR}u:\q y_1...\q y_m[C(y_1,...,y_m) \f D] \f \q y_1...\q y_m[E \f D]$,
$C$ is positive in $E$, $C$ is not free in $D$, and \bf v \rm is not free in $A_j$ $1 \leq j \leq n$.
We may assume that \bf v, \rm $C$ are not free in $B_j$ $1 \leq j \leq n$. By the induction
hypothesis, we have $x_1:B_1,...,x_n:B_n \v_{TTR}u: \q y_1..\q y_m[C(y_1,...,y_m) \f D]\f \q y_1...\q
y_m[E \f D]$, and so $x_1:B_1,...,x_n:B_n \v_{TTR}(Y)u:A$. $\Box$ 
\end{itemize}

\section{$\q$-positive types} 

\subsection{Properties of $\q$-positive types}

\bf Definition \rm We define two sets of types, the set $\O^+$  of
$\q$-positive types, and the set $\O^-$  of $\q$-negative types in the
following way : 
\begin{itemize}
\item[] - If $A$ is an atomic type, then $A \in \O^+$, and $A \in \O^-$ ;  
\item[] - If $T^+ \in \O^+$, and $T^- \in \O^-$, then, $T^- \f T^+ \in \O^+$, and $T^+
\f T^- \in \O^-$ ;  
\item[] - If $T^+ \in \O^+$, then $\q xT^+ \in \O^+$ ;  
\item[] - If $T^+ \in \O^+$, then $\q XT^+ \in \O^+$ ; 
\item[] - If $T^- \in \O^-$, then $\q xT^- \in \O^-$ ; 
\item[] - If $T^- \in \O^-$, and $X$ is not free in $T^-$, then $\q XT^- \in \O^-$- ; 
\item[] - If $T^+ \in \O^+$, $x_1,...,x_n$ first order variables, $t_1,...,t_n$ terms,
$C$ an n-ary predicate symbol which appears and is positive in $T^+$, then
$\m C x_1...x_n T^+ <t_1,...,t_n> \in \O^+$. 
\end{itemize}

\bf Remarks \rm \\
1) A least fixed point is not a $\q$-negative type. \\
2) If $T^+ \in \O^+$, then all least fixed points of $T^+$ are positives. Therefore, by
3)-Theorem 3.2, if $\G \v_{TTR}t:T^+$, then $t$ is normalizable.

\begin{lemma} 
Let $T^-,T'^- \in \O^-$, $T^+,T'^+ \in \O^+$, and $X$ a predicate variable or predicate
symbol.\\ 
1) If $X$ is positive (resp. negative) in $T^-$, then $T^-[T'^-/X] \in \O^-$
(resp. $T^-[T'^+/X] \in \O^-$).\\ 
2) If $X$ is positive (resp. negative) in $T^+$, then $T^+[T'^+/X] \in \O^+$ (resp.
$T^+[T'^-/X] \in \O^+$).\\ 
3) If $T[F/X] \in \O^+$ (resp. $T[F/X] \in \O^-$), then $T \in \O^+$ (resp. $T \in
\O^-$). 
\end{lemma}  
\bf Proof \rm 1), 2) By induction on $T^-$ and $T^+$. \\
3) By induction on $T$. $\Box$ \\

\bf Definition \rm With each type $T$ of $TTR$, we associte the set $Fv_2(T)$ of free
predicate variables and free predicate symbols of $T$.

\begin{theo} Let $T^- \in \O^-$, and $T^+ \in \O^+$.\\
1) If $T^- \inc A$, then $A \in \O^-$, and $Fv_2(A) \inc Fv_2(T^-)$.\\
2) If $B \inc T^+$, then $B \in \O^+$, and $Fv_2(B) \inc Fv_2(T^+)$.
\end{theo}  
\bf Proof \rm We argue by induction on the length of the derivations $T^- \inc A$, and
$B \inc T^+$.  Let us look at the rule used in the last step.\\ 
1) The only non-trivial case is $(\m_d)$. \\
Then $T^-=T'[\m C x_1...x_n T'<y_1,...,y_n>/C(y_1,...,y_n)][t_1/x_1,...,t_n/x_n]$, and \\
$A=\m C x_1...x_n T'<t_1,..,t_n>$. Since $T^- \in \O^-$, then, by Lemma 5.1,
$\m C x_1...x_nT'<y_1,...,y_n> \in \O^-$, which is impossible. \\ 
2) The only non-trivial cases are :
\begin{itemize}
\item[] - $(\m_d)$ : then $B=D[\m C x_1...x_n
D<y_1,...,y_n>/C(y_1,...,y_n)][t_1/x_1,...,t_n/x_n]$, and $T^+=\m C x_1...x_n
D<t_1,...,t_n>$. Since $T \in \O^+$, then $D \in \O^+$, and so, by Lemma 5.1, $B \in
\O^+$, and $Fv_2(B)=Fv_2(D)-\{ C \}=Fv_2(T^+)$.  
\item[] -  $(\m_g)$ : then
$B=\m C x_1...x_n D<t_1,...,t_n>$, $T^+=E[t_1/x_1,...,t_n/x_n]$, and
$D[E/C(x_1,...,x_n)] \inc E$. Since $T^+ \in \O^+$, then $E \in \O^+$, and, by the
induction hypothesis, $D[E/C(x_1,...,x_n)] \in \O^+$, and
$Fv_2(D[E/C(x_1,...,x_n)]) \inc Fv_2(E)$. By Lemma 5.1, we have $D \in \O^+$, and
$Fv_2(D)- \{ C \} \inc Fv_2(D[E/C(x_1,...,x_n)])  \inc Fv_2(E)$, and so $B \in \O^+$,
and $Fv_2(B)= Fv_2(D)- \{ C \} \inc Fv_2(D[E/C(x_1,...,x_n)]) \inc Fv_2(E)=Fv_2(T^+)$.
$\Box$
\end{itemize}

\subsection{ The $TTR_0$ type system}

We define on the types of $TTR$ a binary relation $\inc_0$ by the following way : \\
$A \inc_0$ B if and only if $A \inc B$, and in the proof we use only the weak version of $(\q i_g)$ :
\begin{center}
$(\q i_{g_0})$ $\displaystyle\frac{A[G/v] \inc_0 B} {\q vA \inc_0 B}$ 
\end{center}
where $G$ is a term if $v$ is an individual variable, and $G$ is a predicate variable or
a predicate symbol having the same arity of $v$ if $v$ is a predicate variable.

\begin{lemma} If $A \inc_0 B$, then, for every sequence of terms and/or formulas \bf G
\rm, $A \bf [G/v] \it \inc_0 B \bf [G/v]$, and we use the same proof rules.
\end{lemma}
\bf Proof \rm Same proof as 1)-Lemma 4.2. $\Box$

\begin{lemma} Let $A$ be an arrow type, and $Rep(A)=\q \bf v \it (A_g \f A_d)$.\\
1) If $A \in \O^-$ (resp. $A \in \O^+$), then $A_g \in \O^+$, and $A_d \in \O^-$ (resp.
$A_g \in \O^-$, and $A_d \in \O^+$). \\ 
2) $A \inc_0 Rep(A)$, and $Rep(A) \inc_0 A$.
\end{lemma}
\bf Proof \rm By induction on $A$. $\Box$      
                                                                                           
\begin{lemma} If $T^- \in \O^-$, $T^+ \in \O^+$, and $T^- \inc T^+$, then $T^- \inc_0
T^+$. 
\end{lemma}
\bf Proof \rm By induction on the length of the derivation $T^- \inc T^+$. $\Box$ \\

\bf Definition \rm We denote by $TTR_0$, the $TTR$ type system whithout the rules (5),
(7), (8) and by replacing the rule $(\inc)$ by :
\begin{center}
$( \inc_0 )$ $\displaystyle\frac{\G \v_{TTR_0}t:A \quad A \inc_0 B} {\G\v_{TTR_0}t:B}$ 
\end{center}

\begin{theo} Let $A_1,...,A_n \in \O^-$, $\G=x_1:A_1,...,x_n:A_n$, $A \in \O^+$, and
$t$ a normal $\l$-term. If $\G \v_{TTR}t:A$, then $\G \v_{TTR_0}t:A$, and in this typing
each variable is assigned of a $\q$-negative type, and each $u \in STE(t)$ is typable of
a $\q$-positive type. 
\end{theo} 
\bf Proof \rm We argue by induction on $t$.
\begin{itemize}
\item[]Ê- If $t=x_i$  $1 \leq i \leq n$ , then $\q \bf v$\rm$A_i \inc A$, and \bf v
\rm is not free in $\G$. Since $A_i \in \O^-$, then $\q \bf v$\rm$A_i \in \O^-$, and, by
Lemma 5.4, $\q \bf v$\rm$A_i \inc_0 A$. Therefore $\G \v_{TTR_0}t:A$.   
\item[] - If $t=\l xu$, then $\G,x:B\v_{TTR}u:C$, $\q \bf v$\rm$(BÊ\f C) \inc A$, and
\bf v \rm is not free in $\G$. Since $\q \bf v$\rm$(B \f C)$ is an arrow type, then, by
Lemma 4.5, $A$ is an arrow type. If $Rep(A)=\q \bf v'$\rm$(A_g \f A_d)$, then, by
1)-Lemma 5.3, $A_g \in \O^-$, and $A_d \in \O^+$. By Theorem 4.1, there is a sequence
\bf G, \rm such that $A_g \inc B \bf [G/v]$ \rm, and $C \bf [G/v]$\rm$\inc A_d$. By
2)-Lemma 4.2, we have $\G,x:B \bf [G/v]$\rm$\v_{TTR}u:C \bf [G/c]$ \rm, and, by Lemma
4.9, $\G,x:A_g \v_{TTR}u:A_d$. By the induction hypothesis, we have $\G,x:A_g
\v_{TTR_0}u:A_d$,Êand so, by 2)-Lemma 5.3, $\G \v_{TTR_0}t:A$.   
\item[] - If $t=(x_i)u_1...u_k$  $1 \leq i \leq n$  and $k \not = 0$, then
$\q \bf v_0$\rm$A_i \inc C_1 \f B_1$, $\q \bf v_j$\rm$B_i \inc C_{j+1} \f B_{j+1}$
$1 \leq j \leq k-1$, $\q \bf v_k$\rm$B_k \inc A$ where $ \bf v_0,...,v_k$ \rm are not
free in $\G$, and $\G \v_{TTR}u_j:C_j$ $1 \leq j \leq k$. By Theorems 4.1, 5.1,  and
Lemmas 4.4, 5.4, we have  
\begin{itemize}
\item  $A_i=\q \bf v'_0$\rm$A'_i$, $A'_i=C'_1 \f \q \bf v'_1$\rm$B'_1$,
$B_j=C'_{j+1} \f \q \bf v'_{j+1}$\rm$B'_{j+1}$ $1 \leq j \leq k-1$, $C'_j \in \O^+$, and
$\q \bf v'_j$\rm$B'_j \in \O^-$ $1 \leq j \leq k$. 
\item  $C_j \inc C'_j \bf [G_0/v'_0]...[G_{j-1}/v'_{j-1}]$ \rm,
$\q \bf c'_j$\rm$B'_j \bf [G_0/v'_0]...[G_{j-1}/v'_{j-1}] \inc B_j$ $1 \leq j \leq k$,
and $\q \bf v_k$\rm$\q \bf v'_k$\rm$B'_k \bf [G_0/v'_0]...[G_{k-1}/v'_{k-1}]$\rm$\inc_0
A$. 
\end{itemize}
Since $\G \v_{TTR}u_j:C_j$ $1 \leq j \leq k$, then
$\G \v_{TTR}u_j:C'_j \bf [G_0/v'_0]...[G_{j-1}/v'_{j-1}]$, \rm and, by the induction
hypothesis, $\G \v_{TTR_0}u_j:C'_j \bf [G_0/v'_0]...[G_{j-1}/v'_{j-1}]$. It is easy to
check that $\G \v_{TTR_0}t:B'_k \bf [G_0/v'_0]...[G_{k-1}/v'_{k-1}]$, \rm then \\ 
$\G \v_{TTR_0}t: \q \bf v_k$\rm$\q \bf v'_k$\rm$B'_k \bf [G_0/v'_0]...[G_{k-1}/v'_{k-1}]$, \rm
and $\G \v_{TTR_0}t:A$. $\Box$  
\end{itemize}

\section{G\H{o}del transformation}

\subsection{$\perp$-types of $TTR$}

\bf Definition \rm Let $A$ be a type of $TTR$. We say that $A$ is an $\perp$-type if
and only if $A$ is obtained by the following rules : 
\begin{itemize}
\item[] - $\perp$ is an $\perp$-type.
\item[] - If $A$ is an $\perp$-type, then $B \f A$ is an $\perp$-type for every type
$B$. 
\item[] - If $A$ is an $\perp$-type, then $\q vA$ is an $\perp$-type for every variable
$v$. 
\item[] - If $A$ is an $\perp$-type, $C$ an $n$-ary predicate symbol which appears and
is positive in $A$, $x_1,...,x_n$ first order variables, and $t_1,...,t_n$ terms, then
$\m C x_1...x_n A<t_1,...,t_n>$ is an $\perp$-type. 
\end{itemize}

\begin{lemma} If $A$ is an $\perp$-type, and $A \inc B$, then $B$ is an $\perp$-type. 
\end{lemma}         
\bf Proof \rm By induction on the length of the derivation $A \inc B$. $\Box$
                                                                                                                        
\begin{lemma} Let $t$ be a normal $\l$-term, $A_1,...,A_n \in \O^-$, $A \in \O^+$,
$\perp$ does not appear in the types $A_1,...,A_n,A$, and $B_1,...,B_m$ are
$\perp$-types.  If $\G=x_1:A_1,...,x_n:A_n,y_1:B_1,...,y_m:B_m \v_{TTR}t:A$, then
$x_1:A_1,...,x_n:A_n \v_{TTR}t:A$. 
\end{lemma}
\bf Proof \rm We argue by induction on t. 
\begin{itemize}
\item[] - If $t$ is a variable, then $t=x_i$  $1 \leq i \leq n$  or $t=y_i$ $1 \leq i
\leq m$. 
\begin{itemize}
\item The case $t=x_i$ is trivial. 
\item If $t=y_i$, then $\q \bf v$\rm$B_i \inc A$ and \bf v \rm is not free in $\G$.
Since $B_i$ is an $\perp$-type, then, by Lemma 6.1, $A$ is an $\perp$-type, and $\perp$
appears in $A$. A contradictoire. 
\end{itemize} 
\item[] - If $t=\l x_{n+1}t'$, then $\G,x_{n+1}:A_{n+1}\v_{TTR}t':D$, $\q
\bf v$\rm$(A_{n+1} \f D) \inc A$, \bf v \rm is not free in $\G$. Since $A \in \O^+$,
then, by Theorem 5.1, we have $A_{n+1} \in \O^-$, $D \in \O^+$, and $Fv_2(\q
\bf v$\rm$(A_{n+1} \f D)) \inc Fv_2(A)$. Therefore $\perp$ does not appear in $A_{n+1}$
and $D$. By the induction hypothesis, we have
$x_1:A_1,...,x_n:A_n,x_{n+1}:A_{n+1}\v_{TTR}t:D$, and so
$x_1:A_1,...,x_n:A_n\v_{TTR}t:A$. 
\item[] - If $t=(x)u_1...u_k$  $k \geq 1$, then two case can be see : 
\begin{itemize}  
\item If $x=y_i$ $1 \leq i \leq m$, then, by Corollary 4.1, we have $\q \bf
v_0$\rm$B_i \inc C_1 \f D_1$, $\q \bf v_j$\rm$D_j \inc C_{j+1} \f D_{j+1}$ $1 \leq i
\leq k-1$, $\q \bf v_k$\rm$D_k \inc A$, where $ \bf v_0,...,v_k$ \rm are not free in $A$
and $\G$, and $\G \v_{TTR}u_j:C_j$ $1 \leq j \leq k$. Since $B_i$ is an $\perp$-type,
then, by Lemma 6.1, $D_j$ $1 \leq j \leq k$ and $A$ are $\perp$-types, and $\perp$
appears in $A$. A contradictoire. 
\item If $x=x_i$  $1 \leq i \leq n$, then, by Corollary 4.1,
we have $\q \bf v_0$\rm$A_i \inc C_1 \f D_1$, $\q \bf v_j$\rm$D_j \inc C_{j+1}Ê\f
D_{j+1}$ $1 \leq j \leq k-1$, $\q \bf v_k$\rm$D_k \inc A$, where $ \bf v_0,...,v_k$ \rm
are not free in $A$ and $\G$, and $\G \v_{TTR}u_j:C_j$ $1 \leq j \leq k$. Since $A_i \in
\O^-$, then, by Theorem 5.1, we have $C_j \in \O^+$, $D_i \in \O^-$ $1 \leq j \leq k$,
and $Fv_2(C_j) \bigcup Fv_2(D_j) \inc Fv_2(A_i)$ $1 \leq j \leq k$. Therefore $\perp$
does not appear in $C_j$ $1 \leq j \leq k$. By the inductive hypothesis, we have
$x_1:A_1,...,x_n:A_n\v_{TTR}u_j:C_j$ $1 \leq j \leq k$, and so
$x_1:A_1,...,x_n:A_n\v_{TTR}t:A$. $\Box$  
\end{itemize} 
\end{itemize}
 
\subsection{G\H{o}del transformations} 

\bf Definition \rm With each predicate variable $X$, we associate a finite no empty set
of predicate variables $V_X=\{ X_1,...,X_r \}$ having the same arity of $X$, such that :
if $X \not = Y$, then $V_X \bigcap V_Y= \vi$. With each $n$-ary predicate
variable $X$, and with each sequence of individual variables $x_1,...,x_n$, we assosiate
a formula $F_X$ such that :
\begin{itemize} 
\item[]Ê- $F_X$ is an $\perp$-type ;  
\item[]Ê- $F_X$ does not contain any predicate symbol ; 
\item[]Ê- the free variables of $F_X$ are among $x_1,...,x_n$ and the elements of $V_X$. 
\end{itemize} 
For each formula $A$, we define the formula $A$* by the following induction way : 
\begin{itemize}
\item[]Ê- If $A=C(t_1,...,t_n)$, and $C$ is a predicate symbol, then $A$*=$A$. 
\item[]Ê- If $A=X(t_1,...,t_n)$, and $X$ is a predicate variable, then
$A$*=$F_X[t_1/x_1,...,t_n/x_1]$.  
\item[]Ê- If $A=B \f C$, then $A$*=$B$*$\f C$*. 
\item[]Ê- If $A=\q xB$, then $A$*=$\q xB$*. 
\item[]Ê- If $A=\q XB$, then $A$*=$\q X_1...\q X_r B$*, where $V_X=\{ X_1,...,X_r \}$. 
\item[]Ê- If $A=\m C x_1...x_n D<t_1,...,t_n>$, then $A$*=$\m C x_1...x_n
D$*$<t_1,...,t_n>$.  
\end{itemize} 
$A$* is called the G\H{o}del transformation of $A$.\\

\bf Remark. \rm In order to show that the above transformation is well defined, we need
to prove the following Lemma :

\begin{lemma} Let $C$ be a predicate variable or a predicate symbol, and $A$ a type of
$TTR$. If $C$ is positive in $A$ (resp. negative in $A$), then $C$ is positive in $A$*
(resp. negative in $A$*). 
\end{lemma}
\bf Proof \rm By induction on $A$. $\Box$ 

\begin{lemma} 1) If $A \inc_0 B$, then $A$* $\inc_0 B$*, and we use the same proof
rules.\\ 
2) If $\G \v_{TTR_0}t:A$, then $\G$*$\v_{TTR_0}t:A$*, and we use the same typing
rules. 
\end{lemma}
\bf Proof \rm By induction on the length of the derivation $A \inc_0 B$ (resp. $\G 
\v_{TTR_0}t:A$). $\Box$   

\begin{corollary} Let $D \in \O^+$, and $t$ a normal $\l$-term. If $\v_{TTR}t:D$, then
$\v_{TTR}t:D$*.
\end{corollary}
\bf Proof \rm By induction on the length of the derivation $\v_{TTR}t:D$, and we use
Theorem 5.2 and Lemma 6.4. $\Box$

\section{Storage operators}

\subsection{Definition of storage operators}

\bf Definitions \rm \\ 
1) Let $T$ be a closed $\l$-term, and $D,E$ two closed types of $TTR$
(resp. $TTR^{\di}$). We say that $T$ is a storage operator for the pair of types
$(D,E)$ if and only if for every $\l$-term $t$ with $\v_{TTR}t:D$ (resp.
$\v_{TTR^{\di}}t:D$), there are $\l$-terms $\t_t$ and $\t'_t$ such that $\t_t
\simeq_{\b} \t'_t$, $\v_{TTR}\t'_t:E$ (resp. $\v_{TTR^{\di}}\t'_t:E$), and for every
$\th_t \simeq_{\b} t$, $(T)\th_tf \p (f)\t_t[t_1/x_1,...,t_n/x_n]$, where
$Fv(\t_t)=\{ f,x_1,...,x_n \}$ and $t_1,...,t_n$ are $\l$-terms which depend on
$\th_t$. \\ 
2) If $D=E$, we say that $T$ is a storage operator for the type $D$.\\

\bf Examples \rm The type of recursive integers is the formula : 
\begin{center}
$N^r[x]=\m N x \Ph (N,x)<x>$ 
\end{center}
where 
\begin{center}
$\Ph (N,x)=Ê\q X \{ \q y(Ny \f Xsy),X0 \f Xx \}$
\end{center} 
($s$ is a unary function symbol for successor and $0$ is a constant symbol for zero).\\
For each integer $n$, we define the recursive integer $\sou{n}$ by induction :
$\sou{0}=\l f\l xx$ and $\sou{n+1}=\l f\l x(f)\sou{n}$. Let $\sou{N}$ be the set of recursive
integers.\\  
We have $\sou{N} =\{ t$ / $t$ is a closed normal $\l$-term / $\v_{TTR}t:N^r[s^n(0)]$, $n
\geq 0 \}$ (see [19]).\\ 
Let $\sou{s}=\l n\l f\l x(f)n$. It is easy to check that $\sou{s}$ is a $\l$-term for
successor, and $\v_{TTR}\sou{s}:\q y(N^r[y] \f N^r[sy])$.\\
Define \\
$T_1=(Y)H$ where
$H=\l x\l y((y)\l z(G)(x)z)\d$, $G=\l x\l y(x)\l z(y)(\sou{s})z$, and
$\d=\l f(f)\sou{0}$ ;\\ 
$T_2=\l \n(\n)\r \t \r$ where $\t=\l d\l f(f) \sou{0}$, and $\r=\l y\l z(G)(y)z \t z$,\\
then, for every $\th_n \simeq_{\b} \sou{n}$, $(T_i)\th_nf \p (f)(\sou{s})^n\sou{0}$
$(i=1,2)$. \\
Therefore, for every $n \geq 0$, $T_1$ and $T_2$ are storage operators for
$N^r[s^n(0)]$.\\

\bf Typing of $T_1$ \rm \\

We use in the typing the G\H{o}del transformation with $V_X=\{ X \}$, and
$F_X=\neg X(x_1,...,x_n)$ for every second order variable $X$ of arity $n$.
\begin{itemize} 
\item We have $\v_{TTR}\sou{0}:N^r[0]$, then $\v_{TTR}\d:\neg \neg N^r[0]$.
\item We have $\v_{TTR}\sou{s}:\q y(N^r[y] \f N^r[sy])$, then \\
$x:\neg \neg N^r[y],y:\neg N^r[sy],z:N^r[y] \v_{TTR}(y)(\sou{s})z:\perp$ ; hence :\\
$x:\neg \neg N^r[y],y:\neg N^r[sy] \v_{TTR}(x)\l z(y)(\sou{s})z:\perp$ ; therefore :\\
$\v_{TTR} G:\q y(\neg \neg N^r[y]  \f \neg \neg N^r[sy])$.
\item We have $y:\Ph$*$(N,x) \v_{TTR} y:\q y(Ny  \f \neg \neg N^r[sy]),\neg \neg N^r
[0] \f \neg \neg N^r[x]$ ; thus : \\
$x:\q x(Nx \f \neg \neg N^r[x]),y:\Ph$*$(N,x),z:Ny \v_{TTR}(G)(x)z:\neg \neg N^r[sy]$ ;
therefore : \\
$x: \q x(Nx \f \neg \neg N^r[x]),y:\Ph$*$(N,x) \v_{TTR} \l z(G)(x)z:Ê\q y(Ny \f \neg
\neg N^r[sy])$ ; hence \\ 
$x:\q x(Nx \f \neg \neg N^r[x]) \v_{TTR} \l y((y)\l z(G)(x)z)\d:\q x(\Ph$*$(N,x) \f \neg
\neg N^r[x])$ ; therefore :\\  
$\v_{TTR} H:\q x(Nx \f \neg \neg N^r[x]) \f \q x(\Ph$*$(N,x) \f \neg \neg N^r[x])$.
\end{itemize}
And finally $\v_{TTR}T_1:\q x \{ N^r$*$[x] \f \neg \neg N^r[x] \}$. \\

\bf Typing of $T_2$ \rm \\

We use in the typing the G\H{o}del transformation with $V_X=\{ X,X' \}$, and \\
$F_X=X(x_1,...,x_n),X'(x_1,...,x_n) \f \perp$ for every second order variable $X$ of arity $n$.\\ 
Let $R=\q X \q y \{ (X,X \f \neg \neg N^r[0],X \f \neg \neg N^r[y]),X \f \neg \neg N^r[sy] \}$, $D=R
\f \neg \neg N^r[0]$, and $F[x]=R,D,R \f \neg \neg N^r[x]$. \begin{itemize}
\item $\v_{TTR}\l f(f)\sou{0}:\neg \neg N^r[0]$ ; therefore : $\v_{TTR}\t:X \f \neg
\neg N^r[0]$, and $\v_{TTR}\t:R \f \neg \neg N^r[0]$.
\item By the previous typing, we have $\v_{TTR}G:\q y(\neg \neg N^r[y] \f \neg
\neg N^r[sy])$ ; hence : \\ 
$y:X,X \f \neg \neg N^r[0],X \f \neg \neg N^r[y],z:X \v_{TTR}(G)(y)z \t z:\neg
\neg N^r[sy]$ ; therefore $\v_{TTR}\r:R$.  
\item Check that $\Ph$*$(\l xF[x]/N,x) \inc F[x]$.\\
$\Ph$*$(\l xF[x]/N,x)=$ \\
$\q X \q X' \{ \q y(F[y],Xsy,X'sy \f
\perp),(X0,X'0 \f \perp) \f (Xx,X'x \f \perp) \}$ ; \\ 
therefore by specifying $Xx$ by $R$, and $X'x$ by $\neg N^r[x]$ ; we obtain :\\ 
$\Ph$*$(\l xF[x]/N,x) \inc \q y(F[y],R,\neg N^r[sy] \f \perp),(R,\neg N^r[0] \f \perp) \f (R,\neg
N^r[x] \f \perp)$. We need to check  that $R \inc \q y(F[y],R,\neg N^r[sy] \f \perp)$, this is
absolutely true.  
\end{itemize}  
Therefore $N^r$*$[x] \inc F[x]$ and $\n:N^r$*$[x]\v_{TTR}\n:R,D,R \f
\neg \neg N^r[x]$ ; then :\\  $\n:N^r$*$[x]\v_{TTR}(\n)\r \t \r:\neg \neg N^r[x]$ ; and finally
$\v_{TTR}T_2:\q x \{ N^r$*$[x] \f \neg \neg N^r[x] \}$.

\subsection {General Theorem} 

\begin{theo} Let $D,E$ be two $\q$-positive closed types of $TTR$, such that $\perp$
does not appear in $E$. If $\v_{TTR}T:D$* $\f \neg \neg E$, then $T$ is a storage
operator for the pair $(D,E)$.
\end{theo}
\bf Proof \rm It is a consequence from the following Theorem :

\begin{theo} Let $D,E$ be two $\q$-positive closed types of $TTR^{\di}$, such that
$\perp$ does not appear in $E$. If $\v_{TTR}T:D$* $\f\neg \neg E$, then $T$ is a storage
operator for the pair $(D,E)$.
\end{theo}

Indeed:

\begin{lemma} 1) If $T \in \O^+$ (resp. $T \in \O^-$) then $T^{\di} \in \O^+$ (resp.
$T^{\di} \in \O^-$). \\
2) For each G\H{o}del transformation * of $TTR$, there is a G\H{o}del
transformation *$'$ of $TTR^{\di}$ such that : for every type $D$ of $TTR$,
$D$*$^{\di}=D^{\di}$*$'$.  
\end{lemma}
\bf Proof \rm 1) By induction on $T$.  \\
2) *$'$ is the restiction of * on the types of $TTR^{\di}$.  $\Box$ \\
                                                                            
Let $t$ be a normal $\l$-term, such that $\v_{TTR}t:D$. If $\v_{TTR}T:D$*$\f \neg \neg
E$, then, by Theorem 3.3, $\v_{TTR^{\di}}T:D$*$^{\di} \f \neg \neg E^{\di}$. By 2)-Lemma
7.1, there is a G\H{o}del transformation *$'$, such that
$\v_{TTR^{\di}}T:D^{\di}$*$' \f \neg \neg E^{\di}$. Therefore, there are $\l$-terms
$\t_t$ and $\t'_t$, such that $\t_t \simeq_{\b} \t'_t$, $\v_{TTR^{\di}}\t'_t:E^{\di}$, and
$(T)tf \p (f)\t_t[t_1/x_1,...,t_n/x_n]$. By 2)-Corollary 6.1, we have $\v_{TTR}t:D$*,
then $f:\neg E \v_{TTR}(T)tf:\perp$, and $f:\neg E
\v_{TTR}(f)\t_t[t_1/x_1,...,t_n/x_n]:\perp$. Therefore $f:\neg E
\v_{TTR}(f)\t'_t:\perp$, and, by Corollary 4.1, $\v_{TTR}\t'_t:E$. $\Box$ \\ 

\bf We give the proof of Theorem 7.2 in a particular case. \rm \\

Let $N^r=\m N [\q X \{ N \f X,X \f X \} ]$, and * the G\H{o}del transformation with
$V_X=\{ X \}$, and $F_X= \neg X(x_1,...,x_n)$ for every second order variable $X$ of
arity $n$. \\
We will prove that : If $\v_{TTR^{\di}}T:N^r$*$\f \neg \neg N^r$, then $T$
is a storage operator for $N^r$.\\
Because of : if $t$ is a closed normal $\l$-term with $\v_{TTR^{\di}}t:N^r$, then
$t=\sou{n}$ for a certain integer $n$, and it is suffies to prove that : If
$\v_{TTR^{\di}}T:N^r$*$\f \neg \neg N^r$, then, for every $n \geq 0$, there is an $m
\geq 0$ and $\t \simeq_{\b} \sou{m}$, such that, for every $\l$-term $\th_n \simeq_{\b}
\sou{n}$, there is a substitution $\s$, such that $(T)\th_nf \sim (f)\s(\t)$. 

\begin{lemma} If $\G'=\G,x:N^r$*$\v_{TTR^{\di}}(x)u_1...u_n:\perp$, then $n=3$, and
there is a type $G$, such that $\G'\v_{TTR^{\di}}u_1:N^r$*$\f \neg G$,
$\G'\v_{TTR^{\di}}u_2:\neg G$, and $\G'\v_{TTR^{\di}}u_3:G$. 
\end{lemma}
\bf Proof \rm By Corollary 4.1, we have $\q \bf v_0$\rm$N^r$*$\inc A_1 \f B_1$,
$\q \bf v_i$\rm$B_i \inc A_{i+1} \f B_{i+1}$  $1 \leq i \leq n-1$, $\q \bf v_n$\rm
$B_n \inc \perp$, $\bf v_0,...,v_n$ are not free in $N^r$* and $\G$, and
$\G'\v_{TTR^{\di}}u_i:A_i$  $1 \leq i \leq n$. Since $\q \bf v_0$\rm$N^r$*$\inc A_1 \f
B_1$, then, by Theorem 4.1, there is a formula $F$, such that $A_1 \inc N^r$*$\f \neg
F$ and $\neg F \f \neg F \inc B_1$. We have also $\q \bf v_1$\rm$B_1 \inc A_2 \f B_2$,
then $\q \bf v_1$\rm$(\neg F \f \neg F) \inc A_2 \f B_2$, and, by Theorem 4.1, there is
a sequence of formulas $\bf F_1$, \rm such that $A_2 \inc \neg F\bf [F_1/v_1]$ \rm and
$\neg F\bf [F_1/v_1]$\rm$\inc B_2$. Now, since $\q \bf v_2$\rm$B_2 \inc A_3 \f B_3$, we
have $\q \bf v_2$\rm$(\neg F\bf [F_1/v_1]$\rm$) \inc A_3 \f B_3$, and, by Theorem 4.1,
there is a sequence of formulas $\bf F_2$, such that $A_3 \inc F\bf [F_1/v_1,F_2/v_2]$
\rm and $\perp \inc B_3$. By Corollary 4.1, we have $n=3$. Let
$G=F\bf [F_1/v_1,F_2/v_2]$. \rm Since $\bf v_1,v_2$ are not free in $N^r$* and $\G$,
we deduce $\G'\v_{TTR^{\di}}u_1:N^r$*$\f \neg G$, $\G'\v_{TTR^{\di}}u_2:\neg G$, and
$\G'\v_{TTR^{\di}}u_3:G$. $\Box$ \\   

Let $n \geq 0$.\\

\bf Definition \rm An $n$-special application $\th$ is a function from $\{ 0,1,...,n
\}$ to $\L$ with the following properties : $\th(0) \p \sou{0}$ and
$\th(m+1) \p \l f_m \l x_m(f_m)\th(m)$ $0 \leq m \leq n-1$. Ê 

\begin{lemma} For every $\th_n \simeq_{\b} \sou{n}$, there is an $n$-special application
$\th$, such that $\th(n)=\th_n$.
\end{lemma}
\bf Proof \rm Easy. $\Box$ \\  

\bf Definitions \rm \\ 
1) Let $0 \leq m \leq n$ and
$\bf u$\rm$=u_{m,1},u_{m,2},u_{m,3},...,u_{n-1,1},u_{n-1,2},u_{n-1,3}$ a sequence of $\l$-terms.
We denoted by $x_{m,\bf u}$ \rm a constant which does not appear in $\bf u$.\rm \\ 
2) Let $\th$ be an $n$-special application. The $n$-special substitution $S_{\th}$ is
the function on the set $\L$ defined by induction : 
\begin{itemize}
\item[] - If $u=x$, then $S_{\th}(x)=x$ ; 
\item[] - If $u=\l xv$, then $S_{\th}(u)=\l yS_{\th}(v[y/x])$ where $y \not
\in Fv(\th(n))$ ;  
\item[] - If $u=(v)w$, then $S_{\th}(u)=(S_{\th}(v))S_{\th}(w)$ ;
\item[] - If $u=x_{m,\bf u}$,\rm then \\ 
$S_{\th}(u)=\th(m)[S_{\th}(u_{m,1})/f_m,S_{\th}(u_{m,2})/x_m,...,
S_{\th}(u_{n-1,1})/f_{n-1},S_{\th}(u_{n-1,2})/x_{n-1}]$.
\end{itemize}  
An $n$-special substitution is the application $S_{\th}$ associated to a
some $n$-special application $\th$. 

\begin{lemma} Let $\{ U_i \p V_i \}_{1 \leq i \leq r}$ be a sequence of head reductions
such that : \\
$V_i=(x_{m,\bf u}$\it$)u_1u_2u_3$ $0 \leq m \leq n$,
$[U_{i+1}=(u_1)x_{m-1,u_1,u_2,u_3,\bf u}$\it$u_3$ if $m \not = 0$, and
$U_{i+1}=(u_2)u_3$ if $m=0$], and $S_{\th}$ an $n$-special substitution. For every $1
\leq i \leq r$, $S_{\th}(U_1) \sim S_{\th}(V_i)$.  
\end{lemma} 
\bf Proof \rm We argue by induction on $i$.  \\
The case $i=0$ is a consequence of Theorem 2.1. \\
Assume that is true for $i$, and prove it for $i+1$.  \\
If $V_i=(x_{m,\bf u}$\rm$)u_1u_2u_3$ 0$ \leq m \leq n$, then \\
$S_{\th}(V_i)=(\th(m)[S_{\th}(u_{m,1})/f_m,S_{\th}(u_{m,2}/)x_m,...,S_{\th}(u_{n-1,1})/f_{n-1},
S_{\th}(u_{n-1,2})/x_{n-1}])$\\ 
$S_{\th}(u_1)S_{\th}(u_2)S_{\th}(u_3)$.
\begin{itemize}
\item[] - If $m \not = 0$, then $\th(m) \p \l f_{m-1}\l x_{m-1}(f_{m-1})\th(m-1)$, \\
and $S_{\th}(V_i) \sim
(S_{\th}(u_1))\th(m-1)[S_{\th}(u_{m-1,1})/f_{m-1},S_{\th}(u_{m-1,2})/x_{m-1},...,$\\
$S_{\th}(u_{n-1,1})/f_{n-1},S_{\th}(u_{n-1,2})/x_{n-1}]) S_{\th}(u_3)=S_{\th}(U_{i+1})$. 
\item[] - If $m=0$, then
$\th(m) \p \l f \l xx$, and
$S_{\th}(V_i) \sim (\l f \l xx)S_{\th}(u_1)S_{\th}(u_2)S_{\th}(u_3)
\sim (S_{\th}(u_2))S_{\th}(u_3)=S_{\th}(U_{i+1})$. 
\end{itemize} 
By the induction hypothesis we have $S_{\th}(U_1) \sim S_{\th}(V_i)$, then
$S_{\th}(U_1) \sim S_{\th}(U_{i+1})$, and, by Theorem 2.1, $S_{\th}(U_1)
\sim S_{\th}(V_{i+1})$. $\Box$ \\  

\bf Definition \rm A context $\G=f:\neg N,x_{n,\bf
u_0}$\rm$:N^r$*$,x_{m_1,\bf u_1}$\rm$:N^r$*$,...,x_{m_s,\bf u_s}$\rm$:N^r$* where $0 \leq
m_j \leq n$, $1 \leq j \leq s$, is called $n$-good.

\begin{lemma} There is a sequence of head reductions $\{ U_i \p V_i \}_{1 \leq i \leq
r}$ such that : 
\begin{itemize}
\item[] - $U_1=(T)x_n f$ and $V_r=(f)\t$ where $\t \simeq_{\b} \sou{l}$ for some $l
\geq 0$ ; 
\item[] - $V_i=(x_{m,\bf u}$\rm$)u_1u_2u_3$ $0 \leq m \leq n$, and \\
$U_{i+1}=(u_1)x_{m-1,u_1,u_2,u_3,\bf u}$\rm$u_3$ if $m \not = 0$, and $U_{i+1}=(u_2)u_3$
if $m=0$ ;
\item[] - For every $1 \leq i \leq r$, there is an $n$-good context $\G_i$ such that 
$\G_i\v_{TTR^{\di}}V_i:\perp$.
\end{itemize}
\end{lemma}
\bf Proof \rm Since $\v_{TTR^{\di}}T: N^r$*$\f \neg \neg N^r$, then $x_n:N^r$*$,f:\neg
N^r\v_{TTR^{\di}}(T)x_n f:\perp$, and, by Corollary 4.3 and Lemma 7.2, we have $(T)x_n
f \p V_1$ where $V_1=(f)\t$ or $V_1=(x_n)u_1u_2u_3$. \\
Assume that we have the head reduction $U_k \p V_k$ and $V_k \not = (f)\t$. Then $V_k
=(x_{m,\bf u}$\rm$)u_1u_2u_3$ $0 \leq m \leq n$, and, by the induction hypothesis, there
is an $n$-good context $\G_k$ such that $\G_k \v_{TTR^{\di}}(x_{m,\bf
u}$\rm$)u_1u_2u_3:\perp$. By Lemma 7.2, there is a type $G$, such that $\G_k
\v_{TTR^{\di}}u_1:N^r$*$\f \neg G$,  $\G_k \v_{TTR^{\di}}u_2:\neg G$, and $\G_k 
\v_{TTR^{\di}}u_3:G$.  
\begin{itemize}
\item[] - If $m=0$, let $U_{k+1}=(u_2)u_3$. Let $\G_{k+1}=\G_k$. We have $\G_{k+1}
\v_{TTR^{\di}}U_k:\perp$.  
\item[] - If $m \not = 0$, let $U_{k+1}=(u_1)x_{m-1,u_1,u_2,u_3,\bf u}$\rm$u_3$. The
variable $x_{m-1,u_1,u_2,u_3,\bf u}$ is not used before. Indeed, if it is, by Lemma
7.4, the $\l$-term $(T)\sou{n} f$ is not solvable. That is impossible because $f:\neg
N^r \v_{TTR^{\di}}(T)\sou{n}f:\perp$. Therefore
$\G_{k+1}=\G_k,x_{m-1,u_1,u_2,u_3,\bf u}$\rm$:N^r$* is an $n$-good context and
$\G_{k+1}\v_{TTR^{\di}}U_{k+1}:\perp$.  
\end{itemize}
By Corollary 4.3 and Lemma 7.2, we have $U_{k+1} \p V_{k+1}$ where $V_{k+1}=(f)\t$ or
$V_{k+1}=(x_{s,\bf v}$\rm$)v_1v_2v_3$ $0 \leq s \leq n$.\\ 
This constraction always terminates. Indeed, if not, by Lemma 7.4, the $\l$-term $(T)\sou{n}f$ is not
solvable. That is impossible because $f:\neg N^r\v_{TTR^{\di}}(T)\sou{n}f:\perp$.\\ 
Therefore there is $r \geq 0$ and an $n$-good context $\G_r$ such that $V_r=(f)\t$ and
$\G_r\v_{TTR^{\di}}V_r:\perp$. By Lemma 6.2, we have $\t \simeq_{\b} \sou{l}$ for some $l
\geq 0$.  $\Box$ \\  

Let $\th_n$ be a $\l$-term such that $\th_n \simeq_{\b} \sou{n}$. By Lemma 7.3, let
$\th$ be an $n$-special application such that $\th (n)=\th_n$. Let $S_{\th}$ the
$n$-special substitution associated to $\th$. By Lemma 7.4, we have for every $1 \leq i
\leq r$, $(T)\th_nf \sim S_{\th}(V_i)$. In particular, for $i=n$, $(T)\th_n f \sim
S_{\th}((f)\t)= (f)S_{\th}(\t)$. Then $T$ is a storage operator for $N^r$. $\Box$

\end{document}